\documentclass[12pt]{article}
\newtheorem{thm}{Theorem}
\newtheorem{cor}[thm]{Corollary}
\newtheorem{example}[thm]{Example}
\newtheorem{definition}[thm]{Definition}
\begin{document}

\title{Remarkable and Reversible Prime Number Patterns}
\author{H. J. Weber\\Department of Physics\\
University of Virginia\\Charlottesville, 
VA 22904, U.S.A.}
\maketitle
\begin{abstract}
Prime number multiplet classifications and patterns 
are extended to negative integers. The extension from 
prime numbers to single prime powers is also studied. 
Prime number septets at equal distance are given.
It is also shown that each class of generalized twin 
primes of the classification contains a positive 
fraction of all prime pairs. 
\end{abstract}
\leftline{MSC: 11N05, 11N32, 11N80}
\leftline{Keywords:  Prime number triplets, quintets, 
regular multiplets.} 


\section{Introduction}

The basic classification~\cite{hjw1} of prime number 
twins relies on the following sieve principle. Every 
third number from 1 to $\infty$ is divisible by 3.  
If the greatest common divisor $(D,3)=1,$ then one 
of 3 odd numbers at distance $2D$ is divisible by 3. 
Now 3 can be replaced by any odd prime number. This 
generalized sieve principle then yields exceptional 
prime triplets (and more general prime number 
multiplets) in Theor.~7 of Ref.~\cite{hjw1} or 
Theor.~2.3 of Ref.~\cite{hjw2} and Cor.~10 
of~\cite{hjw1} explaining the often noted empirical 
rule that longer (than 3) prime number sequences at 
equal distance $2D$ must have $3|D$ also (so $D$ is 
excluded from Theor.~7 and Theor.~2.3).            

As in Refs.~\cite{hjw1},\cite{hjw2},\cite{hjw3}, 
we ignore as trivial the prime pairs $(2,p)$ of odd 
distance $p-2$ with $p$ any odd prime. In the 
following prime number multiplets will consist of 
odd primes only. We generally follow the notations 
of Refs.~\cite{hjw1},\cite{hjw2}.  
   
\section{Extension to Negative Integers}

In the context of prime number generating polynomials it 
is often useful to include negative values of the argument 
and the polynomial. This motivates extending the twin prime 
classification to negative integers. We find at once that 
the basic classification of twin primes in Theor.~2 of 
Ref.~\cite{hjw1} or Theor.~2.2 of Ref.~\cite{hjw2}, which 
generalizes to all prime number multiplets, remains valid 
when negative integers are included. The running integer 
variable $a$ then is simply allowed to assume negative 
values as well. 

The classes of odd prime number twins at distance $2D$
\cite{hjw2},  
\begin{eqnarray}\nonumber
I: p_{f,i}&=&2a\pm D,~D=1, 3, 5,\ldots\\\nonumber
II: p_{f,i}&=&3(2a-1)\pm D,~D=2, 4, 8, 10, 14,\ldots\\
III: p_{f,i}&=&2a+1\pm D,~D=6, 12, 18,\ldots,
\end{eqnarray}     
contain all infinitely many prime pairs ($p_i<p_f$).  
Special prime pairs are $3, 3+2D$ with $D$ of class II 
and median $3+D\neq 3(2a-1).$ However, if $D$ is odd, 
then $3+D=2a$ is even, and $2a\pm D$ is a pair in class 
I. All others are special pairs of class II with $D$ even. 
The proof is essentially the same as for Theor.~2 of 
Ref.~\cite{hjw1} and Theor.~2.2 of Ref.~\cite{hjw2}. 

The second classification~\cite{hjw2} of prime pairs in 
terms of arithmetic progressions of conductor $6$ also 
extends to negative numbers. 

When exceptional prime triplets are examined one realizes 
that they sometimes continue as triplets to the left, and 
they reverse by multiplying by $-1$. 

{\bf Example~1.} Exceptional triplets often become 
quintets (at most, composed of a triplet to the left and 
another to the right), such as 
\begin{eqnarray*}
(-7,-5,~3,~11,~13);&\rm{~and~}&(-13,-11,-3,~5,~7);\\
(-13,-5,~3,~11,~19);&\rm{~and~}&(-19,-11,-3,~5,~13);\\
(-17,-7,~3,~13,~23);&\rm{~and~}&(-23,-13,-3,~7,~17);\\
(-37,-17,~3,~23,~43);&\rm{~and~}&(-43,-23,-3,~17,~37).
\end{eqnarray*} 
Exceptional quintets sometimes turn into exceptional 
reversible nonets (bi-quintet or superquintet composed of 
a quintet to the left and another to the right), such as 
\begin{eqnarray}\nonumber
&&(-43,-31,-19,-7,5,17,29,41,53)\\
&&(-53,-41,-29,-17,-5,7,19,31,43) 
\end{eqnarray}  
at equal distance 12.

{\bf Theorem~2.} {\it (i) There is at most one prime 
number quintet with distances $[2D,2D,2D,2D]$ for 
given $D=1, 2, 4, 5,\ldots$ and $(3,D)=1$ composed 
of a triplet starting at $3$ going to the right and 
another to the left. 

(ii) When the distances are $[2d_1,2d_2]$ with 
$3|d_2-d_1$ and\\$3\not|d_1$, there is at most one 
prime number quintet $3-2d_1-2d_2,3-2d_2,3,3+2d_1,
3+2d_1+2d_2$ or $3-2d_1-2d_2,3-2d_1,3,3+2d_2,
3+2d_1+2d_2$ for given natural numbers} $d_1,d_2.$  

The proof is essentially the same as for Theor.~7 
of Ref.~\cite{hjw1} or Theor.~2.3 of Ref.~\cite{hjw2}.

The first line of Example~1 exhibits (ii) for the 
distance pattern $[2,8,8,2],$ while the others are 
at equal distances $8, 10, 20,$ respectively, and 
correspond to (i). Remarkably, there is no case 
where the $[2,8]$ distance pattern moves left to 
become $[2,8,2,8].$   

Exceptional prime number multiplets generalize 
similarly. Now, 3 can be replaced by any odd 
prime number $p>3,$ and $3|D$ is necessary to 
exclude the bi-triplets (quintets) of Theor.~2.
 
{\bf Corollary~3.} {\it For any prime $p>3$ there 
is at most one $(2p-1)-$tuple $p-2(p-1)D,\ldots,
p-2D, p, p+2D,\ldots, p+2(p-1)D$ at a given 
distance} $2D, 3|D, p\not|D.$ 

This bi-$p$-tuple is composed of a $p$-tuple to 
the right and one to the left, both starting at 
$p.$ The sequence just ahead of Theor.~2 is a 
case for $p=5$, and $D=6.$

The proof is essentially the same as for Cor.~10 
of Ref.~\cite{hjw1} and is applied to the left 
and right. 

{\bf Example~4.} At distance $2\cdot 3\cdot 5$ a 
reversible decuplet is 
\begin{eqnarray*}
-157, -127, -97, -67, -37, -7, 23, 53, 83, 113.
\end{eqnarray*}
The reversed 6-tuple $7,37,67,97,127,157$ at 
equal distance 30 has maximum length, being 
different from the exceptional 7-tuple below. 
Going left can at most yield 5 more primes. 
Thus, the decuplet is actually one short of 
optimal. Is there an 11-tuple starting from 7 
at some distance divisible by 6? The exceptional 
septet of G. Lemaire (1909)~\cite{led}
\begin{eqnarray*}
7, 157, 307, 457, 607, 757, 907,
\end{eqnarray*}
is followed by other exceptional septets 
\begin{eqnarray*}
&&7, 2767, 5527, 8287, 11047, 13807, 16567;\\
&&7, 3457, 6907, 10357, 13807, 17257, 20707;\\
&&7, 9157, 18307, 27457, 36607, 45757, 54907;\\
&&7, 14197, 28387, 42577, 56767, 70957, 85147;\\ 
&&7, 21247, 42487, 63727, 84967, 106207, 127447;\\
&&7, 63607, 127207, 190807, 254407, 318007, 381607;\\
&&7, 76717, 153427, 230137, 306847, 383557, 460267;\\
&&7, 117427, 234847, 352267, 469687, 587107, 704527;\\
&&7, 134257, 268507, 402757, 537007, 671257, 805507
\end{eqnarray*}
at larger distances $2760=6\cdot 4\cdot 
5\cdot 23,~3450=6\cdot 5^2\cdot 23,~6\cdot 5^2
\cdot 61,~6\cdot 5\cdot 11\cdot 43, 6\cdot 20\cdot 
3\cdot 59, 6\cdot 10^2\cdot 106, 6\cdot 5\cdot 2557, 
6\cdot 10\cdot 19\cdot 103, 6\cdot 5^3\cdot 179,$ 
respectively. None continues to the left. Are there 
exceptional (reversible) quintets, septets and 
11-, 13-, 17-$\ldots-$tuples at infinitely 
many distances? The distances are all multiples 
of the mandatory $6$; this follows from Theor.~2 
above or Theor.~2.3 of Ref.~\cite{hjw2}. It is 
no accident that distances of exceptional 
septets are divisible by 5. This is needed to 
avoid quintets. This observation generalizes as 
follows.  

{\bf Corollary~5.} {\it An exceptional p-tuple of 
primes starting with the prime p at equal distance 
$D$ must have $p'|D$ for all primes} $p'<p.$ 

{\bf Proof.} This is required to avoid all $p'-$ 
tuples with one number divisible by $p'$ by 
Cor.~10 of Ref.~\cite{hjw1} or Cor.~3. 

{\bf Example~6.} Exceptional quintets at equal 
distance $6k$ are too numerous to be listed. They 
occur for $k=1,2,7,8,16,21,42,71,79,\\99,\ldots.$ 
  
The first exceptional $11-$plet is at equal
distance $1536160080=210\cdot 8\cdot 13\cdot 37\cdot
1901:$
\begin{eqnarray}\nonumber
&&11,~1536160091,~3072320171,~4608480251,~6144640331,
\\\nonumber&&7680800411,~9216960491,~10753120571,~
12289280651,\\&&13825440731,~15361600811.
\end{eqnarray}
The next six (at distances $2\cdot 3^4\cdot 144379,~
2100\cdot 23\cdot 519763,\ldots$) are
\begin{eqnarray}\nonumber
&&11,~4911773591,~9823547171,~14735320751,~19647094331,
\\\nonumber&&24558867911,~29470641491,~34382415071,~
39294188651,\\\nonumber&&44205962231,~49117735811;\\
\nonumber&&11,~25104552911,~50209105811,~75313658711,
~100418211611,\\\nonumber&&125522764511,~150627317411,
~175731870311,~200836423211,\\\nonumber&&225940076111,
~251045529011;
\\\nonumber&&11,~75275138671,~150550279331,~225825418991,
\\\nonumber&&301100558651,~376375698311,~451650837971,
~526925977631,\\\nonumber&&602201117291,~677476256951,
~752751396611;
\\\nonumber&&11,~83516678501,~167033356991,~250550035481,
\\\nonumber&&334066713971,~417583392461,~501100070951,
~584616749441,\\\nonumber&&668133427931,~751650106421,
~835166784911;
\\\nonumber&&11,~100070721671,~200141443331,~300212164991,
\\\nonumber&&400282886651,~500353608311,~600424329971,
~700495051631,\\\nonumber&&800565773291,~900636494951,
~1000707216611;
\\\nonumber&&11,~150365447411,~300730894811,
~451096342211,~601461789611,\\\nonumber&&751827237011,~
902192684411,~10525581131811,~1202923579211,\\&&1353289026611,
~1503654474011;
\end{eqnarray}
The first 13-tuple comes at the enormous distance
$9918821194590$ and is
\begin{eqnarray}\nonumber
&&13,~9918821194603,~19837642389193,~29756463583783,\\
\nonumber&&39675284778373,~49594105972963,~59512927167553,
\\\nonumber&&69431748362143,~79350569556733,~89269390751323,
\\\nonumber&&99188211945913,~109107033140503,~119025854335093.
\end{eqnarray}

Other patterns such as  $[6,6,6,6]: 5,11,17,23,29$ 
are exceptional, but not $[6,6,6]: 61,67,73,79$ 
which is rather common. The distance pattern $[4,2,4]$ 
is very common, too, probably repeats infinitely 
often, appears embedded in $[6,4,2,4,6]: 31,37,41,
43,47,53$ and doubled in $[4,2,4,2,4]: 7,11,13,
17,19,23$ that repeats as $97,101,103,~107,109,
113.$  

Another chaotic property of primes are their gaps. 
The largest gap in the interval $(1,98)$ is $8.$ 
Gaps increase erratically and go way up to $34$ 
in $(1300, 1400)$ and again in $(2110,2200),$ not  
to be reached again until much later in $(8390,8502).$ 
Up to about $10 000,$ rare gap values are large 
single prime gaps, such as $2\cdot 13, 2\cdot 17,$ 
less rare are $2\cdot 5, 2\cdot 7, 2\cdot 11$ and 
common are $2^n\cdot 3,~n=1,2,3,\ldots$ and $2^m
\cdot 3^n\cdot 5.$  
    
Finally, we list a few rules for forming new and 
remarkable patterns of prime number multiplets. 

{\bf Example~7.} From $641,643,647,653,659$ at the  
distances\\$[2,4,6,6]$ we drop 643 forming the 
quartet $641,647,653,659$ at equal distance 6. From 
$601,607,613,617,619$ we drop 617 to get 
$601,607,613,619$ at equal distance $[6,6,6].$ 

{\bf Contraction:} Omitting an intermediate prime 
number leads to a contraction of the multiplet. 
If the distance pattern is\\$[\dots,n_1,n_0,n_2
\ldots]$ then omitting the intermediate prime 
yields\\$[\ldots,n_1,n_0+n_2,\ldots].$ 

{\bf Insertion:} If there exists an intermediate 
prime number it can be re-inserted to yield a 
multiplet with a longer distance pattern: 
$[\ldots,n_1,n_2,\dots]\to [\ldots,n_1,n_0,n_2,
\dots]$ or $[\ldots,n_0,n_1,n_2,\dots]$ or 
$[\ldots,n_1,n_2,n_0,\dots].$ 
 
{\bf Omission of prime at start or end} leads 
to shorter multiplet and distance pattern.  

\section{Admission of Prime Powers}

When single prime powers are admitted as equivalent 
to prime numbers, this amounts to inserting extra 
numbers from which multiplets can start to the left 
and right, with dramatic consequences. Now, there 
are far fewer exceptional multiplets because they 
may repeat starting at a prime power and, within 
positive numbers, they may go left as well.  

{\bf Example~8.} The triplet $3, 7, 11$ at 
distances $[4,4]$ now repeats as $19, 23, 3^3$ 
and goes on to 31, forming a quartet at equal 
distance 4. The exceptional triplet $3, 5, 13$ 
at distances $[2, 8]$ now repeats as a quintet 
$17, 19, 3^3, 29, 37$ with distance pattern 
$[2,8,2,8],$ which is new. The exceptional triplet 
$3, 5, 7$ repeats as a longer quintet $5, 7, 3^2, 
11, 13$ followed, remarkably, by another one 
$23, 5^2, 3^3, 29, 31$ and the triplet 
$79, 3^4,83.$ In the former cases, we have 
triplets to the left and right starting at the 
prime powers $3^2, 3^3,$ respectively. Both 
have distance patterns $[2,2,2,2]$ and are 
almost exceptional, as this pattern never 
repeats. 

{\bf Corollary~9.} {\it The quintet $-13,-5,3,11,19$ 
at equal distance $8$ in Example~1 stays exceptional.  
The quintets at distances $[2,2,2,2]$ in Example~5 
do not repeat further, etc.}  

That's why there are only the following quartets  
$73,3^4,89,97;\\6553,3^8=6561,6569,6577$ at equal 
distance $8.$

{\bf Proof.} For the quintet to repeat requires that 
the two prime powers $3^m, 5^n$ be within 4 times the 
distance 8, or $3^m\approx 5^n$ within $\leq 2^5,$ 
\begin{eqnarray}
\frac{3^m}{5^n}=1+{\cal O}\left(\frac{2^5}{5^n}
\right),~m,n\to\infty, 
\end{eqnarray} 
where the constant in ${\cal O}$ is $\leq 1.$ 
This does not happen for exponents $m,n\leq 9.$ 
In fact, after $5-3=2,$ these prime powers have 
their close encounter 
\begin{eqnarray}
3^3-5^2=2,~\frac{3^{3m}}{5^{2m}}=(1.08)^m\to
\infty,~m\to\infty,
\end{eqnarray}
to diverge exponentially from each other at the 
rate $1.08.~\diamond$

{\bf Lemma~10.} {\it Two prime powers $p_2^{m_2}=
p_1^{m_1}+2D,~D\geq 1,~p_2>5$ have at most one 
close encounter and never meet again.} 

{\bf Proof.} Because 
\begin{eqnarray}
\frac{p_2^{m_2}}{p_1^{m_1}}>1,~\left(\frac{
p_1^{m_1}}{p_2^{m_2}}\right)^n=(1+\varepsilon)^n
\to\infty,~\varepsilon>0,  
\end{eqnarray}
they diverge exponentially from each other.
~$\diamond$ 

{\bf Example~11.} The numbers $9$ and $7$ meet as 
$3^2-7=2$ and then diverge at the rate 
$\frac{3^{2m}}{7^m}=(1.2857\ldots)^m,$ similarly 
$3^4-79=2$ at $\frac{3^{4m}}{79^m}=(1.025\ldots)^m,$ 
$5^3-11^2=4$ at $\frac{5^{3m}}{11^{2m}}=(1.033
\ldots)^m,$ etc. 

{\bf Corollary~12.} {\it Let $p_2^{m_2},p_1^{m_1}
=p_2^{m_2}+2D$ be the only adjacent prime power 
members of a quintet with equal distance pattern 
$[2D,2D,2D,2D].$ If $p_2-p_1>2D$ then this 
quintet is exceptional.}

{\bf Proof.} The sequence of the form $p_1^{m_1+1},
p_2^{m_2+1}=p_1^{m_1+1}+2D,\ldots$ is ruled out 
because this implies $(p_2-p_1)p_2^{m_2}
=2D(p_1+1)<2Dp_2,$ whereas $(p_2-p_1)p_2^{m_2}>2D
p_2^{m_2},$ q.e.a. 
 
A sequence of the form $p_1^{m_1+1},p_3=p_1^{m_1+1}+
2D,p_2^{m_2+1}=p_1^{m_1+1}+4D,\ldots$
is ruled out because this implies $(p_2-p_1)p_2^{m_2}
=2D(p_1+2)<2Dp_2,$ whereas $(p_2-p_1)p_2^{m_2}>2D
p_2^{m_2},$ q.e.a.

A sequence of the form $p_1^{m_1+1},p_3,p_4,p_2^{m_2+1}
=p_1^{m_1+1}+6D,\ldots$ is ruled out because this 
implies $(p_2-p_1)p_2^{m_2}=2D(p_1+3)<2Dp_2,$ whereas 
$(p_2-p_1)p_2^{m_2}>2Dp_2^{m_2},$ q.e.a. 

A sequence of the form $p_1^{m_1+1},p_3,p_4,p_5,
p_2^{m_2+1}=p_1^{m_1+1}+8D,\ldots$ is ruled out 
because this implies $(p_2-p_1)p_2^{m_2}=2D(p_1+4)
\leq 2Dp_2,$ as $p_2\geq p_1+4$ from $p_2-p_1>2D,$ 
whereas $(p_2-p_1)p_2^{m_2}>2Dp_2^{m_2},$ q.e.a.  

Since $m_2\geq 1,~m_1>m_2$ we know that $m_1\geq 2.$ 
So $p_1^{m_1+1}=p_1(p_2^{m_2}+2D)\geq 3p_2^{m_2}+6D.$
If the exponent of $p_2$ stays $m_2$ then the next 
quintet must start like the first: $p_2^{m_2},
p_2^{m_2}+2D=p_1^{m_1},\ldots,$ and then another 
power of $p_1$ cannot occur. So the next quintet 
has to have $p_1^{m_1+1}$ and $p_2^{m_2+1}.$ 
Higher powers are ruled out because they diverge, 
as shown in Lemma~10 and Example~11. This proves 
that no other quintet at equal distance $2D$ 
is possible.~$\diamond$
 
Prime multiplets at such meeting points of prime 
powers are the quintets in Example~8 and Corollary~9, 
or $41,43,47,7^2,53$ at distances $[2,4,2,4]$ that 
is a part-repeater of the longer sequence 
$3,5,7,11,13,17,19,23$, just as $11^2,5^3,127,131$ 
with reversed distance pattern $[4,2,4]$ and 
$163,167,13^2,173$ are repeaters of $7,11,13,\\17$.  
The distance pattern $[4,2,4]$ is very common and 
probably repeats infinitely often, it appears 
also embedded in $[6,4,2,4,6]: 31,37,41,43,47,53$ 
or $[4,2,4,2,4]: 7,11,13,17,19$ that repeats as 
$97,101,103,107,109,113.$  

Another consequence is that some almost-optimal 
prime number generating polynomials become 
optimal, a case in point being 
\begin{eqnarray}
Q_{13}(x)=x^2-3x+43=E_{41}(x-2),~Q_{13}(42)=41^2.       
\end{eqnarray}
A related case is $Q_{14}(j)=$prime, $j=0,\ldots,11$ 
forming a 12-plet, where  
\begin{eqnarray}
Q_{14}(x)=x^2-3x+13,~Q_{14}(12)=11^2, 
\end{eqnarray}
now forming a $13-$tuple, i.e. becoming optimal. 
In addition, $Q_{14}(-x)=x^2+3x+13=$prime for 
$x=0,\ldots, 8;~Q_{14}(-9)=11^2$ forming a 23-plet 
over positive and negative arguments.  
 
\section{Number of Primes in Class I,II,III}

With $a, D$ running in each class I,II,III of 
the classification, it should to be easier to 
settle the question whether or not there 
are infinitely many generalized twins in 
each class than any twin prime conjecture 
with fixed distance $2D.$ 

{\bf Theorem~13.} {\it Each class I,II,III 
contains a positive fraction of all prime 
pairs. Altogether they encompass all prime 
pairs, except for special prime pairs from 
class II, $(3,3+2D)$ with median $3+D\neq 
3(2a-1).$}

{\bf Proof.} We construct the prime pair 
Dirichlet series using\\Golomb's arithmetic 
formula~\cite{sg} for class III 
\begin{eqnarray}
2\Lambda(2a+1-6D)\Lambda(2a+1+6D)=
\sum_{(2a+1)^2-6^2D^2}\mu(d)\log^2d,
\end{eqnarray}
in conjunction with the product theorem for 
Dirichlet series as in Ref.~\cite{hjw3}. 
Let's pick a prime number $p'$ from class III 
and hold it fixed. Write it as $p'=2a+1-6D$ 
and let $a,D$ run so that $p=2a+1+6D$ is a 
matching prime twin in class III for all 
appropriate $a,D$. Now we apply the methods of 
Ref.~\cite{hjw3} to construct the Dirichlet 
series for these pairs of twins. The constraint 
Dirichlet series for this case is  
\begin{eqnarray}\nonumber
q(s)&=&\frac{1}{p'^s}\sum_{a,D}\frac{
\delta_{p',2a+1-6D}}{[2a+1+6D]^s}\\\nonumber
&=&\frac{1}{p'^s}\sum_{D=1}^\infty\frac{1}
{(12D+p')^s}\\&=&\frac{1}{4p'^s}\sum_{\chi_{12}}
\bar\chi_{12}(p')L(s,\chi_{12})-\frac{1}{p'^{2s}},~
\sigma>1, 
\end{eqnarray}
involving an arithmetic progression of conductor 
12. The twin prime series is given by  
\begin{eqnarray}\nonumber
&&2\frac{\log p'}{p'^w}\sum_{D=1}^\infty\frac{
\Lambda(12D+p')}{(12D+p')^w}\\\nonumber&=&-
\frac{\log p'}{2p'^w}\sum_{\chi_{12}}\bar\chi_{12}
(p')\left(\frac{L'}{L}(w,\chi_{12})+\frac{
\chi_{12}(p')\log p'}{p'^w}\right)\\&=&
\lim_{T\to\infty}\frac{1}{2\pi T}\int_{-T}^T 
Z(w-\sigma-it)q(\sigma+it)dt,~\sigma>3.  
\label{a3}
\end{eqnarray}  
Here $Z(s)$ is the twin prime sieve function~\cite{hjw3} 
based on Golomb's formula~\cite{hjw3} 
\begin{eqnarray}
Z(s)=\zeta(s)\frac{d^2}{ds^2}\frac{1}{\zeta(s)}=
\sum_{n=2}^\infty\frac{1}{n^s}\sum_{d|n}\mu(d)\log^2d
\end{eqnarray}
for $\sigma>1.$ Although the product formula~\ref{a3} 
holds for $\sigma>3$ only, the lhs analytically 
continues to the left-hand complex $w-$plane.  
The series on the lhs of Eq.~\ref{a3} contains the 
primes (and powers) of the arithmetic progression 
$12D+p'.$ With $p'=2a+1-6D$ and matching prime 
$2a+1+6D$ the constraint $(2a+1-6D,2a+1+6D$)=1 is 
fulfilled. Here $\chi_{12}$ are the characters 
(mod 12). We now invoke the prime number theorem 
for arithmetic progressions due to Siegel and 
Walfisz which states that, for any conductor 
$q>1,$ the primes are evenly distributed among 
the congruence classes coprime to $q.$ The 
series~(\ref{a3}) has a simple pole at $w=1$ 
with a positive residue. The nonnegative 
coefficients on the lhs allow applying a 
Tauberian theorem implying that the density 
of these twin members in class III is a 
positive fraction of all primes. Letting 
$p'$ run covers all prime pairs of class III.     

We proceed similarly for class II, pick a 
fixed prime $p'\in\{II\},$ write it as 
$3(2a-1)-D$ with $a,D$ running so that 
$3(2a-1)+D=p$ with $3\not|D,~2|D$ is a 
matching twin prime $\in\{II\}.$ The 
constraint series is given by
\begin{eqnarray}
q(s)=\frac{1}{p'^s}\sum_{a=1}^\infty\frac{1}
{[6(2a-1)-p']^s},\sigma>1.
\end{eqnarray} 
The twin prime Dirichlet series for class II 
becomes 
\begin{eqnarray}\nonumber
&&2\frac{\log p'}{p'^w}\sum'_{a>[(p'+6)/12]}
\frac{\Lambda(6(2a-1)-p')}{[6(2a-1)-p']^w}\\
\nonumber&=&-\frac{\log p'}{p'^w}\bigg[\sum_{
\chi_6}\frac{\bar\chi_6(-p')}{2}\frac{L'}
{L}(w,\chi_6)-\sum_{\chi_{12}}\frac{
\bar\chi_{12}(-p')}{4}\frac{L'}{L}(w,\chi_{12})
\bigg]\\&=&\lim_{T\to\infty}\frac{1}{2T}
\int_{-T}^T Z(w-\sigma-it)q(\sigma+it)dt,~\sigma>3,
\end{eqnarray}
where $\chi_6,\chi_{12}$ are the characters 
(mod 6) and (mod 12), respectively, and $[x]$ 
is the largest integer below $x.$. Again, 
the series is over the primes (and their powers) 
of the relevant arithmetic progression. The 
application of the prime number theorem in 
conjunction with the Tauberian theorem leads 
to the same conclusion. 

For class I, $p'=2a-D, 2\not|D, p=2a+D,$ 
Golomb's identity is given by
\begin{eqnarray}
2\Lambda(2a-D)\Lambda(2a+D)=\sum_{d|4a^2-D^2}
\mu(d)\log^2 d.
\end{eqnarray}
The constraint Dirichlet series is defined as 
\begin{eqnarray}\nonumber
q(s)&=&\frac{1}{p'^s}\sum_{a>[p'/4]}^\infty\frac{1}
{(4a-p')^s}\\&=&\frac{1}{2p'^s}\sum_{\chi_4}
\bar\chi_4(-p')[L(s,\chi_4)-\chi_4(1)],
~\sigma>1. 
\end{eqnarray}
The twin prime Dirichlet series becomes 
\begin{eqnarray}\nonumber
&&2\frac{\log p'}{p'^w}\sum_{a>[p'/4]}^\infty
\frac{\Lambda(4a-p')}{(4a-p')^w}\\\nonumber&=&-
\frac{\log p'}{p'^w}\sum_{\chi_4}\bar\chi_4(-p')
\frac{L'}{L}(w,\chi_4)\\&=&
\lim_{T\to\infty}\frac{1}{2T}\int_{-T}^T 
Z(w-\sigma-it)q(\sigma+it)dt,~\sigma>3,
\end{eqnarray}
where $\chi_4$ are the characters (mod 4). 
Again, it has a simple pole at $w=1,$ the 
Tauberian theorem applies and gives the 
corresponding result.$~\diamond$

Next, let us apply this method to the classes 
of the second generalized twin prime 
classification, Theor.~2.3 of Ref.~\cite{hjw2}.  

{\bf Corollary~14.} {\it The subset of twin primes 
$p',p\equiv 1\pmod{6}$ for $a\equiv 0\pmod{3}$ in 
class III is a positive fraction of all prime pairs.}

{\bf Proof.} We have 
\begin{eqnarray}
a=3\alpha,~p'=6\alpha+1-6D,~p=6\alpha+1+6D,
\end{eqnarray}
where $p'$ is held fixed again, while $\alpha$ 
and $D$ run. The constraint Dirichlet series is 
\begin{eqnarray}
q(s)&=&\frac{1}{p'^s}\sum_{D=1}^\infty\frac{1}
{(12D+p')^s}\\&=&\frac{1}{4p'^s}\sum_{\chi_{12}}
\bar\chi_{12}(p')L(s,\chi_{12})-\frac{1}{p'^{2s}},
~\sigma>1.
\end{eqnarray} 
The associated twin prime series is 
\begin{eqnarray}\nonumber
&&2\frac{\log p'}{p'^w}\sum_{D=1}^\infty\frac{
\Lambda(12D+p')}{(12D+p')^w}\\\nonumber&=&-
\frac{\log p'}{2p'^w}\sum_{\chi_{12}}\bar
\chi_{12}(p')\left(\frac{L'}{L}(w,\chi_{12})
+\frac{\chi_{12}(p')\log p'}{p'^w}\right)\\&=&
\lim_{T\to\infty}\frac{1}{2T}\int_{-T}^T 
Z(w-\sigma-it)q(\sigma+it)dt,~\sigma>3,
\end{eqnarray}
where $\chi_{12}$ are the characters (mod 12). 
Again, it has a simple pole at $w=1,$ the 
Tauberian theorem applies and gives the 
corresponding result.$~\diamond$ 

The method applies similarly to other 
classes of the second generalized twin 
prime classification in Theor.~2.3 of 
Ref.~\cite{hjw2}.


\begin{thebibliography}{0}  

\bibitem{hjw1}  H. J. Weber, 2011, ``Exceptional 
prime number twins, triplets and multiplets,'' 
arXiv:1102.3075 [math.NT]. 

\bibitem{hjw2}  H. J. Weber, 2011, ``Regularities 
of twin, triplet and multiplet prime numbers,''  
Global J. of Pure and Applied Math.{\bf 8} (2012), 
arXiv:1103.0447 [math.NT].   

\bibitem{hjw3} H. J. Weber, 2010, ``Generalized 
Twin Prime Formulas,'' Global J. of Pure and 
Applied Math. 6(1), 101-116, arXiv:1106.1054 [math.NT].   

\bibitem{sg} S. W. Golomb, 1970, ``The Lambda 
Method in Prime Number Theory,'' J. Number 
Theory, {\bf 2}, 193-198. 

\bibitem{led} L. E. Dickson, 2005, ``History of 
the Theory of Numbers,'' Vol. I, Dover, Mineola, 
N.Y., p. 426.

\end{thebibliography}
\end{document}